\documentclass[conference]{IEEEtran}

\ifCLASSINFOpdf
\else
\fi
\usepackage{graphics} 
\usepackage{epsfig} 
\usepackage{amsmath} 
\usepackage{amssymb}  
\usepackage[latin1]{inputenc}
\usepackage{amsfonts}
\usepackage{subfigure}

\newtheorem{lemme}{Lemma}
\newtheorem{prop}{Proposition}

\begin{document}

\title{
Non-asymptotic fractional order differentiators via an algebraic parametric method
}

\author{\IEEEauthorblockN{Da-Yan Liu}
\IEEEauthorblockA{Mathematical and \\ Computer Sciences and \\ Engineering Division,\\
King Abdullah University \\of Science and Technology (KAUST), KSA\\
Email: dayan.liu@kaust.edu.sa}
\and
\IEEEauthorblockN{Olivier Gibaru}
\IEEEauthorblockA{ LSIS  (CNRS, UMR 7296),\\ Arts et Métiers
ParisTech \\Centre de Lille,\\
L'\'{E}quipe Projet Non-A,\\
INRIA Lille-Nord Europe, France\\
Email: olivier.gibaru@ensam.eu}
\and
\IEEEauthorblockN{Wilfrid
Perruquetti}
\IEEEauthorblockA{LAGIS (CNRS, UMR 8146),\\ \'{E}cole Centrale de Lille,\\
L'\'{E}quipe Projet Non-A,\\
INRIA Lille-Nord Europe, France\\
Email: wilfrid.perruquetti@inria.fr}}

\maketitle

\begin{abstract}
Recently, Mboup, Join and Fliess \cite{Mboup2007,Mboup2009a} introduced non-asymptotic integer order differentiators by using an algebraic parametric estimation method \cite{FLI03b_mfhsr,FLI08a_garnier}.
In this paper, in order to obtain  non-asymptotic
fractional order differentiators  we apply this
algebraic parametric method to  truncated expansions of fractional Taylor series based on the Jumarie's modified Riemann-Liouville derivative \cite{Jumarie2009a}. Exact and simple formulae for these differentiators are given where a sliding integration window of
a noisy signal involving Jacobi polynomials  is used without complex mathematical deduction.
 The  efficiency
and the stability with respect to
corrupting noises of the proposed fractional order differentiators are shown  in numerical
simulations.
\end{abstract}

\IEEEpeerreviewmaketitle
\section{INTRODUCTION}

Fractional models arise in many practical situations (\cite{Oustaloup10,Oustaloup02} for example). Such fractional order systems may also be used for control purposes: CRONE control
is known to have good robustness properties (see \cite{Oustaloup95, Oustaloup98, Oustaloup99, Oustaloup02}). In order to implement such controller one needs to have a good digital fractional order differentiator from noisy signals, which is the scope of this paper.

The fractional derivative has a long history and has often appeared in
science, engineering and finance (see, e.g., \cite{Oldham,Podlubny2,Fliess98,Fliess97}).
Different from classical integer order derivative, 
there are several kinds of
definitions for the fractional derivative which are generally not equivalent with each
other \cite{Podlubny2,Kilbas}. Among these definitions, the Riemann-Liouville derivative and
the Caputo derivative are often used \cite{Podlubny2}.
Recently, a new modified Riemann-Liouville derivative is
proposed by Jumarie \cite{Jumarie2009a}. This new definition of fractional derivative has two main advantages: firstly comparing with the  Caputo derivative, the function to be differentiated is not necessarily  differentiable, secondly different
from the Riemann-Liouville derivative the Jumarie's modified Riemann-Liouville derivative of a constant is defined to zero. Moreover, a fractional  Taylor series expansion \cite{Jumarie} was invented by using this new definition.
Thanks to these merits, the Jumarie's modified Riemann-Liouville derivative was successfully applied (see, e.g., \cite{Jumarie2008,Jumarie2009b,Wu}). 
The fractional order differentiator is concerned with estimating the
fractional order derivatives of an unknown signal from its noisy observed
data. Because of its importance, various methods have been developed during
the last years (\cite{Chen,Machado} for example). However, the obtained fractional order differentiators were usually based on
 the Riemann-Liouville derivative and
the Caputo derivative. To our knowledge, there is no one based on the Jumarie's modified Riemann-Liouville derivative.

Recent algebraic parametric estimation method for linear systems
\cite{FLI03b_mfhsr,FLI08a_garnier,Tian2008} has been extended to various problems in
signal processing (see, e.g., \cite{FLI03a_mexico,FLI04a_compression_cras,Mboup2009b,Neves2006,Neves2007,Ushirobira2011,Trapero2007,Liu2011d,Liu2008}).
Let us emphasize that this method is algebraic and non-asymptotic, which provides explicit formulae and finite-time estimates. Moreover, it exhibits good
robustness properties with respect to corrupting noises, without the
need of knowing their statistical properties (see \cite{ans,shannon}
for more theoretical details). The robustness properties have
already been confirmed by numerous computer simulations and several
laboratory experiments.
Very recently, this  method was used to solve the ill-posed numerical differentiation problem
in \cite{Mboup2007,Mboup2009a}. The main idea is to
apply a differentiator called integral-annihilator  to a truncated Taylor series expansion which is a local approximation of the signal to be differentiated.
Stable non-asymptotic differentiators of integer order were exactly given by the integrals of noisy
signals involving the Jacobi orthogonal polynomials.
The associated estimation errors  were studied in
\cite{Liu2009,Liu2011b,Liu2011c}.
An extension of these differentiators for multivariate numerical differentiation were proposed in
\cite{Riachy,Riachy2008,Riachy2010}.
However, this method has not been used to estimate
fractional order derivatives.

The aim of this paper is to
introduce non-asymptotic fractional order differentiators for the Jumarie's modified Riemann-Liouville derivative
by using the algebraic parametric method.
 In Section \ref{section1}, we recall how to apply the algebraic parametric method to obtain an integer order differentiator from a truncated Taylor series expansion. Section \ref{section2} begins with the definition of the Jumarie's modified Riemann-Liouville derivative. Then, two fractional order differentiators are obtained by applying integral-annihilators  to  truncated fractional Taylor series expansions with different truncated orders. Moreover, it is shown that the differentiator obtained from higher order truncated expansion can be expressed as an affine combination of the ones obtained from lower order truncated expansion.
Numerical
tests are given in Section \ref{section3}. They help us to show the efficiency
and the stability of the proposed fractional order differentiators.
 Finally, we give some conclusions and
perspectives for our future work in Section \ref{section4}.

\section{METHODOLOGY}\label{section1}

Let $y=x+\varpi$ be a noisy signal observed in an open
interval $I \subset\mathbb{R}$, where
$x \in \mathcal{C}^n(I)$ with $n \in \mathbb{N}$ and the noise\footnote{More generally, the noise is a
stochastic process, which is bounded with certain probability and integrable
in the sense of convergence in mean square (see \cite{Liu2011b}).} $\varpi$ is
bounded and integrable. In this section, we are going to recall how to use an algebraic parametric
method to estimate the integer order derivatives of $x$  \cite{Mboup2007, Mboup2009a}.

For any $t_0
\in I$, we introduce the set $D_{t_0}:=\{t \in \mathbb{R}_+;\, t_0 + t
\in I\}$. By using the famous  Taylor's formula formulated by Hardy
(\cite{Hardy} p. 293), we obtain that $\forall \,t_0 \in I, \ t \in  D_{t_0}, $
\begin{align}
  x(t_0+t)&= \sum_{j=0}^{n}
\frac{ t^j}{j!} x^{(j)}(t_0)+ \mathcal{O}(t^n), \  \text{ as } t
\rightarrow 0.
\end{align}

In a similar  way to classical numerical differentiation methods, the
algebraic parametric method uses  the $n^{th}$ order derivative of a
polynomial to estimate the one of $x$. Here, we consider the
following truncated Taylor series expansion of $x$ on $\mathbb{R}_+$
\begin{align}\label{Eq_Taylor_series_truncated}
\forall\,  t \in \mathbb{R}_+,\ x_n(t_0+t):= \sum_{j=0}^{n}
\frac{t^j}{j!} x^{(j)}(t_0).
\end{align}
Different from the other methods, the $n^{th}$ order derivative of
the approximation polynomial $x_n$ is calculated by applying   algebraic
manipulations to $x_n$ in the operational domain. Precisely, we
apply a differential operator with the following form \cite{Liu2011b}
\begin{align}\label{Eq_annihilator_minimal}
\Pi^{n}_{k,\mu}=\frac{1}{s^{n+1+\mu}} \cdot
\frac{d^{n+k}}{ds^{n+k}}\cdot {s^n},
\end{align}
where $-1 <\mu \in
\mathbb{R}$ and  $k \in \mathbb{N}$.
This  differential operator was introduced in \cite{Mboup2009a} with $\mu \in \mathbb{N}$.

Let us mention that the operation $\frac{d^{n+k}}{ds^{n+k}}\cdot {s^n}$ is used to
annihilate all the terms containing $x^{(j)}(t_0)$ with $0 \leq j
\leq n-1$ in $x_n$ in the operational domain, and the rational term
$\frac{1}{s^{n+1+\mu}}$ permits to obtain a Riemann-Liouville
integral \cite{Loverro} expression  for $x^{(n)}(t_0)$ in the time
domain (see the next section for more details). Thus, the
differential operator $\Pi^{n}_{k,\mu}$ is called
\emph{integral-annihilator for $x^{(n)}(t_0)$ via $x_{n}$} \cite{Mboup2009a}.
Finally, if we replace $x_n$ by $y$ in the obtained integral
expression, then we can obtain the following estimator for
$x^{(n)}(t_0)$ (see \cite{Liu2011b} and \cite{Mboup2009a}  for more
details)
\begin{align}
\begin{split} \label{Eq_minimal_estimator}
& \forall\, t_0 \in I, \ \tilde{y}^{(n)}_{{t_0}}(k,\mu, T):= \\ & \quad \ \ \frac{n!}{
T^{n}} \gamma_{n,k,\mu} \int_{0}^{1}
w_{\mu,k}(\tau) P^{(\mu,k)}_{n}(\tau)\, y(t_0+ T\tau) \,d\tau,
\end{split}
\end{align}
where  $T \in D_{t_0}$, $\gamma_{n,k,\mu}=\frac{1}{\mathrm{B}(n+k+1,n+\mu+1)}$, $\mathrm{B}(\cdot,\cdot)$ is the classical
beta function (\cite{Abramowitz}, p. 258), $P^{(\mu,k)}_{n}$ is the
$n^{th}$ order Jacobi polynomial (\cite{Abramowitz} p. 775) defined
on $[0,1]$  as follows
\begin{equation}\label{Eq_jacobi_poly}
{P}_{n}^{(\mu,k)}(\tau)=\displaystyle\sum_{j=0}^{n}\binom{n+\mu}%
{j}\binom{n+k}{n-j}\left(  {\tau-1}\right)  ^{n-j}\left( \tau\right)
^{j},
\end{equation}
and  $w_{\mu,k}(\tau)=(1-\tau)^{\mu} \tau^{k}$ is the associated
weight function.
Hence, this differentiator depends on three parameters $k$, $\mu$, $T$. Moreover, it is
a non-asymptotic pointwise differentiator using the sliding integration window
$[t_0,t_0+T]$. Since $\tilde{y}^{(n)}_{{t_0}}(k,\mu, T)$ is
obtained by taking an $n^{th}$ order polynomial where $n$ is the
order of the derivative estimated, we call it \emph{minimal Jacobi
differentiator}.

\section{NON-ASYMPTOTIC FRACTIONAL ORDER DIFFERENTIATOR} \label{section2}
In this section, we are going to  propose two  fractional order differentiators
by applying the algebraic parametric
method to different truncated expansions of the fractional  Taylor series based on the Jumarie's modified Riemann-Liouville derivative  \cite{Jumarie2009a,Jumarie}.

\subsection{The Jumarie's modified Riemann-Liouville derivative}

Let $f$ be a continuous function defined on $\mathbb{R}$, then the Jumarie's modified Riemann-Liouville derivative of $f$ is defined as follows \cite{Jumarie2009a}
\begin{equation}
f^{(\alpha)}(t):=\frac{1}{\Gamma(l-\alpha)}\frac{d^{l}}{dt^{l}}%
\int_{0}^{t}\frac{f(\tau)-f(0)}{\left(t-\tau\right)^{\alpha+1-l}}
\,d\tau,
\end{equation}
where $0\leq l-1 \leq \alpha<l$ with $l\in\mathbb{N}^{\ast}$. This fractional order derivative is in fact defined through the
fractional difference \cite{Jumarie2009a}
\begin{equation}
f^{(\alpha)}(t):=\lim_{h\rightarrow 0} \frac{\Delta^{\alpha}\left[f(t)-f(0)\right]}{h^{\alpha}},
\end{equation}
where $h>0$, $FW\, f(t)=f(t+h)$ and
\begin{equation}
\begin{split}
\Delta^{\alpha}f(t)&=\left(FW-1\right)^{\alpha}f(t)\\&=\sum^{\infty}_{i=0} (-1)^i \binom{\alpha}{i}f\left[t+(\alpha-i)h\right].
\end{split}
\end{equation}
Let us recall that the Riemann-Liouville
derivative is defined as follows
(\cite{Podlubny2} p. 62)
\begin{equation}
{}_{0}{D}_{t}^{\alpha}f(t):=\frac{1}{\Gamma(l-\alpha)}\frac{d^{l}}{dt^{l}}%
\int_{0}^{t}\frac{f(\tau)}{\left(t-\tau\right)^{\alpha+1-l}}
\,d\tau,
\end{equation}
where $0\leq l-1 \leq \alpha<l$ with $l\in\mathbb{N}^{\ast}$.
If we take
$f(t)=t^{n}$ with $n\in\mathbb{N}$ and $t\in\mathbb{R}_{+}^{\ast}$, then we
obtain (see \cite{Podlubny2} p. 72)
\begin{equation}
{}_{0}{D}_{t}^{\alpha}t^{n}=\frac{\Gamma(n+1)}{\Gamma(n+1-\alpha
)}\,t^{n-\alpha},\ \text{ with }\alpha\in\mathbb{R}_{+}.
\end{equation}
Consequently, the Jumarie's modified Riemann-Liouville derivative can be expressed as follows
\begin{equation}\label{Eq_modi_RL}
\forall\, t\in\mathbb{R}_{+}^{\ast},\    f^{(\alpha)}(t)={}_{0}{D}_{t}^{\alpha}f(t)-\frac{t^{-\alpha}}{\Gamma(1-\alpha)}\, f(0).
\end{equation}

One of some useful properties of the Jumarie's modified Riemann-Liouville derivative is the fractional Leibniz derivative rule  \cite{Jumarie2009a}
 \begin{equation}\label{Eq_fractional_product}
 \left(f\, g\right)^{(\alpha)}= f^{(\alpha)}g+g^{(\alpha)}f.
\end{equation}

Moreover, since $x \in \mathcal{C}^n(I)$,  a generalized Taylor series
expansion based on the Jumarie's modified Riemann-Liouville derivative of $x$ is given as follows \cite{Jumarie2009a,Jumarie}:
$\forall \, t_0 \in I$,  $t \in D_{t_0}$,
\begin{equation}\label{Eq_Taylor_series_fractional}
x(t_0+ t)=\sum_{j=0}^n
\frac{t^j}{j!}x^{(j)}(t_0)+ \sum_{j=1}^{+\infty}\frac{
t^{j\gamma+n}}{\Gamma(j\gamma+n+1)}x^{(j\gamma+n)}(t_0),
\end{equation}
where  $\gamma=\alpha-n$ with $n < \alpha \leq n+1$.

Let us remark that the fractional order derivative $x^{(j\gamma+n)}$ in (\ref{Eq_Taylor_series_fractional})
should be understood as $(x^{(n)})^{(j\gamma)}$ which is different from $(x^{(j\gamma)})^{(n)}$ \cite{Jumarie2009a}.
\subsection{Minimal Jacobi fractional differentiator}

By using (\ref{Eq_Taylor_series_fractional}), we take the following truncated fractional Taylor series expansion of $x$ on $\mathbb{R}^+$:
$\forall \, t_0 \in I$,  $t \in \mathbb{R}^+$,
\begin{equation}\label{Eq_Taylor_series_fractional_truncated}
 x_{\alpha}(t_0+
t):=\sum_{j=0}^n \frac{ t^j}{j!}x^{(j)}(t_0)+ \frac{
t^{\alpha}}{\Gamma(\alpha+1)}x^{(\alpha)}(t_0).
\end{equation}
Then, by using the algebraic parametric method we can give the
following proposition.

\begin{prop} \label{Prop_estimator_fractional}
Let $y=x+\varpi$ be a noisy signal observed on an open
interval $I \subset\mathbb{R}$, where
$x \in \mathcal{C}^n(I)$ with $n < \alpha \leq n+1$, $n \in \mathbb{N}$, and $\varpi$ be a
bounded and integrable noise. Then an  estimator for the $\alpha^{th}$ order
derivative value $x^{(\alpha)}(t_0)$ is given by: $\forall \, t_0 \in I$,
\begin{equation}
\begin{split} \label{Eq_estimator_fractional}
 &\tilde{y}^{(\alpha)}_{{t_0}}(k,\mu, T):=\\&
\frac{(n+1)!}{ T^{\alpha}} \gamma_{n,k,\mu,\alpha}
 \int_{0}^{1}
w_{\mu,k}(\tau) P^{(\mu,k)}_{n+1}(\tau)\, y(t_0+  T\tau) \,d\tau,
\end{split}
\end{equation}
where  $T \in D_{t_0}$, $\gamma_{n,k,\mu,\alpha}=\frac{\Gamma(\alpha-n)}{\mathrm{B}(\alpha+1+k,n+\mu+2)}$, and $P^{(\mu,k)}_{n+1}$ is the Jacobi polynomial
defined by (\ref{Eq_jacobi_poly}) with $k \in \mathbb{N}$ and $-1<
\mu \in \mathbb{R}$.
\end{prop}

\noindent{\textbf{Proof.}} By applying the Laplace transform to
$(\ref{Eq_Taylor_series_fractional_truncated})$, we get
\begin{equation}\label{Eq_Taylor_series_fractional_laplace}
\hat{x}_{\alpha}(s)=\sum_{j=0}^n  s^{-(j+1)} x^{(j)}(t_0)+
s^{-(\alpha+1)} x^{(\alpha)}(t_0),
\end{equation}
where $\hat{x}_{\alpha}(s)$ is the Laplace transform of
$x_{\alpha}(t_0+ t)$, and $s$ is the Laplace variable.

We are going to apply some algebraic manipulations to
(\ref{Eq_Taylor_series_fractional_laplace}). Firstly, we apply the
operation $\frac{d^{n+1+k}}{ds^{n+1+k}}\cdot s^{n+1}$ to
(\ref{Eq_Taylor_series_fractional_laplace}) so as to annihilate the
terms containing $x^{(j)}(t_0)$ with $0 \leq j \leq n$ in
$\hat{x}_{\alpha}$. Thus, we get
\begin{align*}
\begin{split}
&\frac{d^{n+1+k}}{ds^{n+1+k}} s^{n+1} \left(\sum_{j=0}^n s^{-(j+1)}
x^{(j)}(t_0)+ s^{-(\alpha+1)} x^{(\alpha)}(t_0) \right)\\
=&(-1)^{n+1+k} \frac{
\Gamma(\alpha+1+k)}{\Gamma(\alpha-n)}\frac{1}{s^{\alpha+k+1}}
x^{(\alpha)}(t_0).
\end{split}
\end{align*}

Secondly, if we apply the Leibniz derivative rule to
$\frac{d^{n+1+k}}{ds^{n+1+k}} s^{n+1}\hat{x}_{\alpha}(s)$, then the
highest order of $s$ in the obtained sum is $n+1$. Hence, we choose the
rational term  $\frac{1}{s^{n+2+\mu}}$ with $\mu >-1$ such that
$\frac{1}{s^{n+2+\mu}}\cdot s^{n+1}=\frac{1}{s^{1+\mu}}$ with
$\mu+1>0$. This permits to obtain a Riemann-Liouville integral \cite{Loverro}
expression of $x_{\alpha}$ in the time domain. Consequently, we  construct an
integral-annihilator for $x^{(\alpha)}(t_0)$ via
$x_{\alpha}$ which corresponds to the operator $\Pi^{n+1}_{k,\mu}$ defined by (\ref{Eq_annihilator_minimal}).

Consequently, by applying the inverse Laplace transform  and the
classical rules of operational calculus, we get
\begin{equation}\label{Eq_1}
\begin{split}
&\mathcal{L}^{-1}\left\{\Pi^{n+1}_{k,\mu} \left(\sum_{j=0}^n
s^{-(j+1)} x^{(j)}(t_0)+
s^{-(\alpha+1)} x^{(\alpha)}(t_0) \right)\right\}(T)=\\
& \quad \quad \quad (-1)^{n+1+k}
\frac{\Gamma(\alpha+1+k)}{\Gamma(\alpha-n)}
\frac{T^{n+\alpha+2+k+\mu}}{\Gamma(n+\alpha+k+\mu+3)}
x^{(\alpha)}(t_0),\\
 &\mathcal{L}^{-1}\left\{\Pi^{n+1}_{k,\mu}
\hat{x}_{\alpha}(s) \right\}(T)=\\
& \quad  \frac{(-1)^{n+1+k}}{\Gamma(n+\mu+2)}
\int_{0}^{T} (T-\tau)^{n+\mu+1} \tau^{n+k+1} x^{(n+1)}_{\alpha}(t_0+
\tau) d\tau.
\end{split}
\end{equation}

By applying a change of variable $\tau \rightarrow T\tau$ and $n+1$
times integrations by parts, we get
\begin{equation}
\begin{split} \label{Eq_estimator_fractional0}
& \frac{(-1)^{(n+1)}T^{\alpha}}{\gamma_{n,k,\mu,\alpha}}\, x^{(\alpha)}(t_0)=\\&
\int_{0}^{1}
\frac{d^{n+1}}{d\tau^{n+1}}\left\{(1-\tau)^{n+\mu+1}
\tau^{n+k+1}\right\} x_{\alpha}(t_0+T\tau)\, d\tau.
\end{split}
\end{equation}
Finally, this proof can be completed by substituting $x_{\alpha}$ in
(\ref{Eq_estimator_fractional0}) by $y$ and applying the Rodrigues
formula (\cite{Abramowitz} p. 785) to the right side of
(\ref{Eq_estimator_fractional0}). \hfill$\Box$

Since the differentiator $\tilde{y}^{(\alpha)}_{{t_0}}(k,\mu, T)$ involves a Jacobi polynomial, and it is obtained by taking the
minimal order ($j=1$) truncated fractional Taylor series expansion
in (\ref{Eq_Taylor_series_fractional}), we call it \emph{minimal Jacobi
fractional differentiator}. The corresponding truncated error part
comes from the truncated term
$\displaystyle\sum_{j=2}^{+\infty}\frac{
t^{(j\gamma+n)}}{\Gamma(j\gamma+n+1)}x^{(j\gamma+n)}(t_0)$. If we take $\alpha=n+1$ in
(\ref{Eq_estimator_fractional}), then it is easy to obtain that
\begin{equation}\label{}
\tilde{y}^{(\alpha)}_{{t_0}}(k,\mu,T)=
\tilde{y}^{(n+1)}_{{t_0}}(k,\mu,T),
\end{equation}
where $\tilde{y}^{(n+1)}_{{t_0}}(k,\mu,T)$ is the minimal Jacobi
differentiator for $x^{(n+1)}(t_0)$ given in
(\ref{Eq_minimal_estimator}). Moreover, in a similar way to the minimal
Jacobi differentiator method, the minimal Jacobi fractional differentiator
$\tilde{y}^{(\alpha)}_{{t_0}}(k,\mu,T)$ can also be obtained by
using the classical orthogonal properties of the Jacobi polynomials to (\ref{Eq_Taylor_series_fractional_truncated})
(see \cite{Liu2011a} for more details). Consequently, as done in
\cite{Liu2011b}, the parameter $k$ which is defined on $\mathbb{N}$
can also be extended to $]-1, +\infty[$.

Finally, let us recall that the affine Jacobi differentiator being an affine combination of the minimal Jacobi differentiators was introduced in
 \cite{Mboup2007, Mboup2009a} by applying the algebraic parametric method to a higher order truncated Taylor series expansion than (\ref{Eq_Taylor_series_truncated}). Hence, the convergence rate for the Jacobi differentiator was improved. Similarly to the integer derivative case, we are going to introduce an Affine Jacobi fractional differentiator in the next subsection.


\subsection{Affine Jacobi fractional differentiator}

In this subsection, we take a new truncated fractional Taylor series
expansion in (\ref{Eq_Taylor_series_fractional}) up to $j=2$: $\forall \, t_0 \in I$,  $t \in \mathbb{R}^+$,
\begin{align}\label{Eq_Taylor_series_fractional_truncated2}
\begin{split}
& x_{2\alpha-n}(t_0+
t)= \sum_{j=0}^n \frac{ t^j}{j!}x^{(j)}(t_0)\\ &  + \frac{
t^{\alpha}}{\Gamma(\alpha+1)}x^{(\alpha)}(t_0)+\frac{
t^{2\alpha-n}}{\Gamma(2\alpha-n+1)}x^{(2\alpha-n)}(t_0).
\end{split}
\end{align}
In order to obtain an estimator for
$x^{(\alpha)}(t_0)$ from this truncated expansion, we introduce the following
differential operator
\begin{align}\label{Eq_annihilator_affine}
\Pi^{n+1}_{k,\mu,\alpha}=\frac{1}{s^{2\alpha+k+\mu+3}} \cdot
\frac{d}{ds}\cdot s^{2\alpha-n+1+k} \cdot
\frac{d^{n+k+1}}{ds^{n+k+1}}\cdot s^{n+1},
\end{align}
where $-1 <\mu \in \mathbb{R}$,  $k \in \mathbb{N}$ and $n<\alpha \leq
n+1$ with $n \in \mathbb{N}$. Then, we give the following proposition.

\begin{prop}
The differential operator $\Pi^{n+1}_{k,\mu,\alpha}$ defined in
(\ref{Eq_annihilator_affine}) is an integral-annihilator for
$x^{(\alpha)}(t_0)$ via $x_{2\alpha-n}$ defined in
(\ref{Eq_Taylor_series_fractional_truncated2}).
\end{prop}

\noindent{\textbf{Proof.}} Applying the Laplace transform to
$(\ref{Eq_Taylor_series_fractional_truncated2})$, we get
\begin{align}\label{Eq_Taylor_series_fractional_laplace2}
\begin{split}
&\hat{x}_{2\alpha-n}(s)= \sum_{j=0}^n  s^{-(j+1)} x^{(j)}(t_0)\\& \quad \quad +
s^{-(\alpha+1)} x^{(\alpha)}(t_0) + s^{-(2\alpha-n+1)}
x^{(2\alpha-n)}(t_0),
\end{split}
\end{align}
where $\hat{x}_{2\alpha-n}(s)$ is the Laplace transform of
$x_{2\alpha-n}(t_0+ t)$.

We are going to apply $\Pi^{n+1}_{k,\mu,\alpha}$ to
(\ref{Eq_Taylor_series_fractional_laplace2}). Firstly, we apply the
operation $\frac{d^{n+1+k}}{ds^{n+1+k}}\cdot s^{n+1}$ to
(\ref{Eq_Taylor_series_fractional_laplace2}) so as to annihilate the
terms containing $x^{(j)}(t_0)$ with $0 \leq j \leq n$ in
$\hat{x}_{2\alpha-n}$. Thus, we obtain
\begin{equation}
\begin{split}\label{Eq_Taylor_series_fractional_laplace3}
&\frac{d^{n+1+k}}{ds^{n+1+k}} s^{n+1} \left(\hat{x}_{2\alpha-n}(s) \right)= \\& \quad (-1)^{n+1+k} \frac{
\Gamma(\alpha+1+k)}{\Gamma(\alpha-n)}\frac{1}{s^{\alpha+k+1}}
x^{(\alpha)}(t_0)\\+&(-1)^{n+1+k} \frac{
\Gamma(2\alpha-n+1+k)}{\Gamma(2\alpha-2n)}\frac{x^{(2\alpha-n)}(t_0)}{s^{2\alpha-n+k+1}}.
\end{split}
\end{equation}

Secondly, we apply the operator $\frac{d}{ds}\cdot
s^{2\alpha-n+1+k}$ to (\ref{Eq_Taylor_series_fractional_laplace3})
so as to annihilate the term containing $x^{(2\alpha-n)}(t_0)$ in
$\hat{x}_{2\alpha-n}$. This yields
\begin{equation}
\begin{split}\label{Eq_Taylor_series_fractional_laplace4}
&\frac{d}{ds}s^{2\alpha-n+1+k}\frac{d^{n+1+k}}{ds^{n+1+k}} s^{n+1}
\left(\hat{x}_{2\alpha-n}(s) \right)=
\\&\quad  (-1)^{n+1+k} \frac{
\Gamma(\alpha+1+k)}{\Gamma(\alpha-n)}(\alpha-n)s^{\alpha-n-1}
x^{(\alpha)}(t_0).
\end{split}
\end{equation}

Thirdly, if we apply the Leibniz  derivative rule to
$\frac{d}{ds}s^{2\alpha-n+1+k}\frac{d^{n+1+k}}{ds^{n+1+k}}
s^{n+1}\hat{x}_{2\alpha-n}(s)$, then the highest order of $s$ in the
obtained sum is $2\alpha+n+k+2$. Hence, we apply the rational term
$\frac{1}{s^{2\alpha+n+k+3+\mu}}$ with $\mu >-1$ such that
$\frac{1}{s^{2\alpha+n+k+3+\mu}}\cdot
s^{2\alpha+n+k+2}=\frac{1}{s^{1+\mu}}$ with $\mu+1>0$. This allows us
to obtain a Riemann-Liouville integral expression of $x_{2\alpha-n}$
in the time domain. Consequently, $\Pi^{n+1}_{k,\mu,\alpha}$ is an
integral-annihilator for $x^{(\alpha)}(t_0)$ via $x_{2\alpha-n}$.
Moreover, we have
\begin{equation}
\begin{split} \label{Eq_2}
&(-1)^{n+1+k} \, \mathcal{L}^{-1}\left\{\Pi^{n+1}_{k,\mu, \alpha} \hat{x}_{2\alpha-n}(s) \right\}(T)=\\
&
\frac{\Gamma(\alpha+1+k)}{\Gamma(\alpha-n)}
\frac{(\alpha-n)\, T^{n+\alpha+3+k+\mu}}{\Gamma(n+\alpha+k+\mu+4)}
x^{(\alpha)}(t_0),
\end{split}
\end{equation}
where $T \in D_{t_0}$. \hfill$\Box$

Similarly to Proposition \ref{Prop_estimator_fractional}, by
substituting $x_{2\alpha-n}$ by $y$ in the obtained
Riemann-Liouville integral $\mathcal{L}^{-1}\left\{\Pi^{n+1}_{k,\mu,
\alpha} \hat{x}_{2\alpha-n} \right\}(T)$, we obtain a new estimator
for $x^{(\alpha)}(t_0)$. We denote it by
$\tilde{y}^{(\alpha)}_{{t_0}}(k,\mu,T,\alpha)$. Then, we get the
following proposition.

\begin{prop} \label{Prop_estimator_fractional_affine}
Let $y$ be a noisy signal defined as in Proposition \ref{Prop_estimator_fractional},
then we have the following affine relation
\begin{align}\label{Eq_estimator_relation}
\begin{split}
\tilde{y}^{(\alpha)}_{{t_0}}(k,\mu,T,\alpha)&=\lambda_{\alpha,k,n}\,
\tilde{y}^{(\alpha)}_{{t_0}}(k,\mu+1,T)\\&  + (1-\lambda_{\alpha,k,n})\,
\tilde{y}^{(\alpha)}_{{t_0}}(k+1,\mu,T),
\end{split}
\end{align}
where $\lambda_{\alpha,k,n}=\frac{2\alpha-n+1+k}{\alpha-n}$.\ The differentiator
$\tilde{y}^{(\alpha)}_{{t_0}}(k,\mu,T,\alpha)$ is thus called \emph{affine Jacobi
fractional differentiator}.
\end{prop}

In order to prove the above proposition, we need the following
lemma.
\begin{lemme} \label{Lemme}
Let $\hat{f}$ be the Laplace transform of an analytic function $f$
defined on $I$, $\Pi^{n+1}_{k,\mu,\alpha}$ and $\Pi^{n+1}_{k,\mu}$
be the integral-annihilators defined in
(\ref{Eq_annihilator_affine}) and (\ref{Eq_annihilator_minimal})
respectively, then we have
\begin{equation}\label{Eq_annihilator_relation}
\begin{split}
&\Pi^{n+1}_{k,\mu,\alpha}\hat{f}(s)=\\
&\quad \quad (2\alpha-n+1+k)\,
\Pi^{n+1}_{k,\mu+1}\hat{f}(s) + \Pi^{n+1}_{k+1,\mu}\hat{f}(s).
\end{split}
\end{equation}
\end{lemme}

\noindent{\textbf{Proof.}} By applying the Leibniz derivative rule, we get
\begin{align*}
&\Pi^{n+1}_{k,\mu,\alpha}\hat{f}(s)\\ =&\frac{1}{s^{2\alpha+k+\mu+3}}
\left( (2\alpha-n+1+k) s^{2\alpha-n+k} \frac{d^{n+k+1}}{ds^{n+k+1}}
s^{n+1}\hat{f}(s)\right) \\ &+ \frac{1}{s^{2\alpha+k+\mu+3}} \left( s^{2\alpha-n+k+1} \frac{d^{n+k+2}}{ds^{n+k+2}}
s^{n+1}\hat{f}(s)\right)\\
=&(2\alpha-n+1+k) \Pi^{n+1}_{k,\mu+1}\hat{f}(s) +
\Pi^{n+1}_{k+1,\mu}\hat{f}(s).
\end{align*}
\hfill$\Box$

\noindent{\textbf{Proof of Proposition
\ref{Prop_estimator_fractional_affine}.}} According to (\ref{Eq_2}),
we obtain
\begin{equation}
\begin{split} \label{Eq_3}
&\mathcal{L}^{-1}\left\{\Pi^{n+1}_{k,\mu, \alpha} \hat{x}_{2\alpha-n}(s)\right\}(T)\\= &\mathcal{L}^{-1}\left\{\Pi^{n+1}_{k,\mu,
\alpha} \left(\sum_{j=0}^n s^{-(j+1)} x^{(j)}(t_0)+
s^{-(\alpha+1)} x^{(\alpha)}(t_0) \right)\right\}(T)\\
= &\mathcal{L}^{-1}\left\{\Pi^{n+1}_{k,\mu, \alpha}
\hat{x}_{\alpha}(s)\right\}(T)\\
=&  (-1)^{n+1+k} \frac{\Gamma(\alpha+1+k)}{\Gamma(\alpha-n)}
 \frac{(\alpha-n) T^{n+\alpha+3+k+\mu}}{\Gamma(n+\alpha+k+\mu+4)}
x^{(\alpha)}(t_0).
\end{split}
\end{equation}
By substituting $x_{\alpha}$ by $x$ in the Riemann-Liouville
integrals $\mathcal{L}^{-1}\left\{\Pi^{n+1}_{k,\mu, \alpha}
\hat{x}_{\alpha}(s)\right\}(T)$ and
$\mathcal{L}^{-1}\left\{\Pi^{n+1}_{k,\mu}
\hat{x}_{\alpha}(s)\right\}(T)$ given in (\ref{Eq_3}) and
(\ref{Eq_1}) respectively, we get
\begin{equation}
\begin{split} \label{Eq_4}
\tilde{x}^{(\alpha)}_{{t_0}}(k,\mu,T,\alpha):= &
\frac{c_{n,k,\mu,\alpha,T} }{\alpha-n}
\mathcal{L}^{-1}\left\{\Pi^{n+1}_{k,\mu, \alpha}
\hat{x}(s)\right\}(T)\\
\tilde{x}^{(\alpha)}_{{t_0}}(k,\mu+1,T):=& c_{n,k,\mu,\alpha,T} \mathcal{L}^{-1}\left\{\Pi^{n+1}_{k,\mu+1}
\hat{x}(s) \right\}(T)
\\
\tilde{x}^{(\alpha)}_{{t_0}}(k+1,\mu,T):=&
-\frac{c_{n,k,\mu,\alpha,T}}{\alpha+1+k}
\mathcal{L}^{-1}\left\{\Pi^{n+1}_{k+1,\mu} \hat{x}(s) \right\}(T),
\end{split}
\end{equation}
where $c_{n,k,\mu,\alpha,T}=(-1)^{n+1+k}
\frac{\Gamma(\alpha-n)}{\Gamma(\alpha+1+k)}
\frac{\Gamma(n+\alpha+k+\mu+4)}{T^{n+\alpha+3+k+\mu}}$.
Hence, by using Lemma \ref{Lemme} we obtain
\begin{align}\label{Eq_estimator_relation2}
\begin{split}
\tilde{x}^{(\alpha)}_{{t_0}}(k,\mu,T,\alpha)& =\lambda_{\alpha,k,n}\,
\tilde{x}^{(\alpha)}_{{t_0}}(k,\mu+1,T)\\& + (1-\lambda_{\alpha,k,n})\,
\tilde{x}^{(\alpha)}_{{t_0}}(k+1,\mu,T).
\end{split}
\end{align}
Finally, this proof can be completed by substituting $x$ by $y$ in $\tilde{x}^{(\alpha)}_{{t_0}}(k,\mu,T,\alpha)$, $\tilde{x}^{(\alpha)}_{{t_0}}(k+1,\mu,T,\alpha)$ and $\tilde{x}^{(\alpha)}_{{t_0}}(k,\mu+1,T,\alpha)$.

\hfill$\Box$


\begin{figure}[h!]
\centering {\includegraphics[scale=0.5]{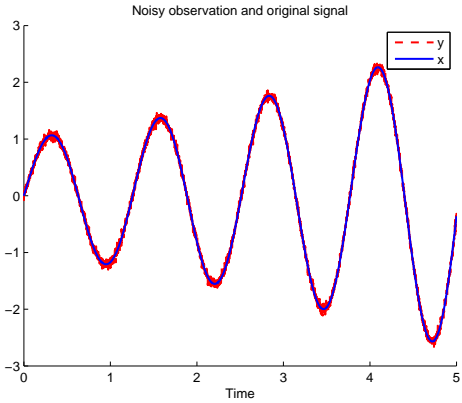}}
\caption{Signal $x$ and
noisy signal $y$.}%
\label{fig_signal}%
\end{figure}

\begin{figure}[h!]
 \centering
 {\includegraphics[scale=0.5]{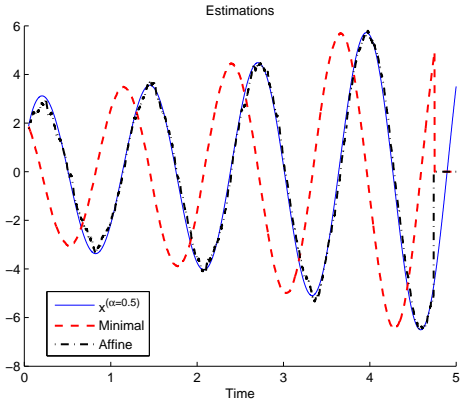}}
\caption{
$\tilde{y}^{(0.5)}_{{t_0}}(0,0,T)$ with $T=0.25$ and
$\tilde{y}^{(0.5)}_{{t_0}}(0,0,T,0.5)$ with $T=0.26$.}
\label{fig_estimations}
\end{figure}

\begin{figure}[h!]
 \centering
 {\includegraphics[scale=0.5]{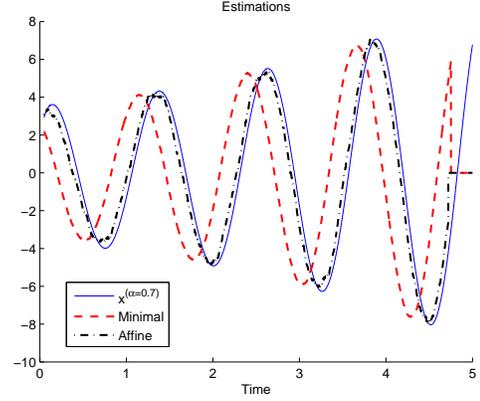}}
\caption{
$\tilde{y}^{(0.7)}_{{t_0}}(0,0,T)$ with $T=0.25$ and
$\tilde{y}^{(0.7)}_{{t_0}}(0,0,T,0.7)$ with $T=0.28$.}
\label{fig_estimations2}
\end{figure}


\section{NUMERICAL SIMULATIONS} \label{section3}
In order to show the efficiency and the stability of the proposed fractional order differentiators,
we give some numerical results in this section.

We assume that $y(t_i)=x(t_i)+\varpi(t_i)$ is the discrete noisy
observation of $x$ on $I=[0,4]$  where $x(t_i)= \exp(0.2t_i)\sin(5t_i)$, $\varpi$
is a zero-mean white Gaussian noise, and $t_i=T_s i$ for
$i=0,\cdots,4 \times 10^3$ with $T_{s}=\frac{1}{10^3}$. The variance of
$\varpi$ is adjusted in such a way that the signal-to-noise ratio
$SNR=10\log_{10}\left(\frac{\sum|y(t_i)|^{2}}{\sum|\varpi(t_i)|^{2}}\right)$
is equal to $SNR=28.07\text{dB}$.  We can see the original signal $x$ and its noisy observation $y$ in Figure \ref{fig_signal}. The exact Jumarie's modified Riemann-Liouville derivative of $x$ can be calculated by using (\ref{Eq_fractional_product}), (\ref{Eq_modi_RL}) and the Riemann-Liouville derivative of the functions $\sin(\cdot)$ and $\exp(\cdot)$ (see \cite{Miller}, p. 83).   We  estimate the
derivatives $x^{(\alpha)}$ with $\alpha=0.5$ and $0.7$ respectively.
We apply the trapezoidal numerical integration method to
approximate the integrals in our differentiators where we use $m+1$ ($T=m T_s$)
discrete observation values in each sliding integration window.

The formal derivatives and their estimated values  are shown in Figure
\ref{fig_estimations} and Figure \ref{fig_estimations2}, where the
parameters $k$ and $\mu$ are set to zero for the Jacobi  fractional
differentiators. On one hand, the error due to the noise for the minimal Jacobi fractional differentiator can be negligible with respect to the truncated error part which produces a time-shift in the estimate.
On the other hand, the noise error for the affine Jacobi fractional differentiator is larger than the one for the minimal Jacobi fractional differentiator, but the truncated error part is much smaller. Hence,  the time-shift is significant reduced.


\section{CONCLUSION} \label{section4}

In this paper, two non-asymptotic
fractional order differentiators called minimal and affine Jacobi fractional differentiators are proposed by applying an
algebraic parametric method to truncated expansions of fractional Taylor series. They can be used to estimate the recently invented  Jumarie's modified Riemann-Liouville derivative. Numerical simulations are given so as to show their
efficiency and stability with respect to corrupting noises.
It is shown that these differentiators can also be obtained by
using the classical orthogonal properties of the Jacobi polynomials \cite{Liu2011a}. By using this method,
a generalized affine Jacobi fractional differentiator will be obtained from a truncated fractional Taylor series expansion with an arbitrary truncated order in the future work. Moreover,
in a similar way to \cite{Liu2011a} and \cite{Liu2011b}, the noise
errors and the truncated errors in these differentiators  will be analyzed. In particular, we will study the influence of
the parameters $k$, $\mu$ and $T$ to these errors so as to give a
guideline  for choosing the optimal parameters.


\end{document}